\documentstyle[11pt]{article}
\renewcommand{\^}{\hspace*{-3pt}{}^}
\def\st{\scriptstyle}

\def\ds{\displaystyle}

\def\nz{\normalsize}
\def\Ce{{\rm C}\hspace*{-1.92mm}\rule{0.5pt}{2.35mm}\hspace*{2mm}}
\def\Re{\,{\rm R}\hspace*{-3.4mm}{\rm I} }
\newcommand{\be}{\begin{equation}\label}
\newcommand{\ee}{\end{equation}}
\newcommand{\ba}{\begin{eqnarray}}
\newcommand{\ea}{\end{eqnarray}}
\newcommand{\ar}{\begin{array}}
\newcommand{\er}{\end{array}}

\newcommand{\up}[2]{\stackrel{\scriptstyle#1}{#2}}
\def\wt{\widetilde}        \def\rar{\rightarrow}
\def\wh{\widehat}          
\def\red{\hspace*{6mm}}    \def\ov{\overline}
\def\p{{\vec{\phantom{p}}\hspace*{-5pt}p}}

\def\a{\alpha}           \def\C{{\cal C}}     \def\k{\kappa}
\def\b{\beta}            \def\F{{\cal F}\,}
\def\g{\gamma}           \def\G{{\cal G}}
\def\s{\sigma}           \def\J{{\cal J}}
\def\D{\Delta}           
\def\d{\delta}           \def\U{{\cal U}}    
\def\la{\lambda}         \def\H{{\cal H}}
\def\eps{\epsilon}       \def\M{{\cal M}}
\def\om{\omega}          
\def\ve{\varepsilon}     \def\R{{\cal R}}
\def\o{\otimes}          \def\P{{\cal P}}

\def\CG{{\Bigl\{\ar{c}\!\st\!I\;\;\,J\;\;K\\ [-5pt] \st m\;m'\,m'' 
\er\!\Bigr\}}}

\def\nn{\nonumber}

\begin{document}
\hyphenation{ope-ra-tor ope-ra-tors ge-ne-ra-tor ge-ne-ra-tors
ge-ne-ra-ting ope-ra-tion ope-ra-tions pre-print quan-tum} 

\thispagestyle{empty}
September~~1996
\hfill{ \begin{tabular}{l} \tt Berlin Sfb288~No.215 \\
 \tt q-alg/9609007 \end{tabular} }

\vspace*{2.5cm}  \begin{center}
{\Large \bf Fusion of $q$-tensor operators : \\ [2mm]
 quasi-Hopf-algebraic point of view.}  
   \end{center}
\vspace{8mm}

\begin{center}  {\sc Andrei G. Bytsko} $^*$   \\ [4mm]
  Institut f\"ur Theoretische Physik, Freie Universit\"at Berlin \\
Arnimallee 14, 14195 Berlin, Germany \\
and \\
Steklov Mathematical Institute, Fontanka 27,\\
 St.Petersburg~~191011, Russia
\end{center}

\vspace*{2.5cm}

\begin{center} {\bf Abstract  } \\ [3mm]

\parbox{13.5cm}{

Tensor operators associated with a given quantum Lie algebra 
$U_q(\J)$ admit a natural description in the $R$-matrix language.
Here we employ the $R$-matrix approach to discuss the problem of 
fusion of tensor operators.  
The most interesting case is provided by the quantum WZNW model, 
where, by construction, we deal with sets of linearly independent 
tensor operators. 
In this case the fusion problem is equivalent to construction of an 
analogue $\F(p)$ of the twisting element $\F$ which is employed in 
Drinfeld's description of quasi-Hopf algebras. 
We discuss the construction of the twisting element $\F(p)$ in a 
general situation and give illustrating calculations for the case 
of the fundamental representation of $U_q(sl(2))$.
} \end{center}

\vfill

\begin{flushleft}   \rule{7 cm}{0.05 mm} \\
$^{*}$ \ { e-mail : \ 
bytsko@physik.fu-berlin.de , \ bytsko@pdmi.ras.ru  }
\end{flushleft}

\newpage
\setcounter{page}{1}

\section{\nz\bf INTRODUCTION}

\subsection{Motivations and notations}

\red The theory of tensor operators has arisen originally as a result of a
group-theoretical treatment of quantum mechanics \cite{WBL}. 
And conversely, the
further development of the representation theory was inspired by 
the physical interpretation of its mathematical content.
Relatively recent appearance of the theory of quantum groups \cite{Dr1} 
has led to development of the theory of $q$-deformed tensor operators 
\cite{BTR,Bie}. 
The latter turned out to be not purely mathematical construction; 
it is employed, in particular, in the description of the quantum  WZNW 
model \cite{AF,CG,BS}.
 
In the present paper we discuss some aspects of the 
{\it fusion} procedure for (deformed) tensor operators in its $R$-matrix 
formulation \cite{Byt}. We consider the special case of the fusion 
scheme -- a construction of a set of basic tensor operators for given irrep 
$\rho^K$ if we are given those for two other irreps $\rho^I$ and $\rho^J$ 
(and $\rho^K$ appears in decomposition of $\rho^I\o\rho^J$).
It turns out that this problem is closely related to Drinfeld's 
construction of quasi-Hopf algebras \cite{Dr2}.
Our aim is to obtain exact prescriptions applicable in practice, but 
to formulate the problem precisely we need to give first rather 
detailed introduction to the subject.

We suppose that the reader is familiar with the notion of Hopf algebra. 
The latter is an associative algebra $\G$ equipped with unit $e\in\G$,
a homomorphism $\D:\G\mapsto\G\o\G$ (the {\em co-product}), an 
anti-automorphism $S:\G\mapsto\G$ (the {\em antipode}) and 
a one-dimensional representation $\eps:\G\mapsto\Ce$ (the {\em co-unit}) 
which obey a certain set of axioms \cite{AS}. 
A quasi-triangular Hopf algebra \cite{Dr1} possesses in addition
an invertible element $R\in\G\o\G$ (the {\em universal $R$-matrix}) obeying
certain relations which, in particular, imply the Yang-Baxter equation.

Throughout the paper we shall use so-called $R$-matrix formalism 
\cite{FRT,KS}. Let us recall that its main ingredients are 
operator-valued matrices ($L$-operators)
\be{0.1}
 L^I_+=(\rho^I\o id)\,R_+ \ , \;\;\;\;   L^I_-=(\rho^I\o id)\,R_-
\ee
and numerical matrices ($R$-matrices)
\be{0.2}
R^{IJ}_+=(\rho^I\o \rho^J)\,R_+ \ , \;\;\;\;
R^{IJ}_-=(\rho^I\o \rho^J)\,R_- \ ,
\ee
with $\rho^I$, $\rho^J$ being irreps of $\G$ and $R_+=R$, 
$R_-=(R^{\,\prime})^{-1}$.
Here and further on $'$ stands for permutation in $\G\o\G$.

Our consideration is restricted to the case of $\G=U_q(\J)$
with $|q|=1$ and with $\J$ being a semi-simple Lie algebra (the case of 
$\G$ being an arbitrary semi-simple modular Hopf algebra needs some 
additional technique; see the discussion in \cite{BS}). For simplicity 
we assume also that $q$ is not a root of unity.

We perform all explicit computations only in
the case of $U_q(sl(2))$, but they can be certainly repeated 
for, at least, $U_q(sl(n))$.
Let us underline also that, although we deal with deformed tensor 
operators and keep index $q$ in some formulae, the classical 
(i.e., non-deformed) theory is recovered  in 
the limit $q=1$ and, therefore, it does not need special comments.

\subsection{(Deformed) tensor operators, generating matrices }

\red Let the given quasi-triangular Hopf algebra $\G$ be a symmetry 
algebra for some physical model. 
This means that the operators corresponding to the
physical variables in this model are classified with respect to their
transformation properties under the adjoint action of $\G$. Recall 
that if $\cal H$ is a certain Hilbert space such that 
$\G\subset{\rm End}\,{\cal H}$, then ($q$-deformed) adjoint action of  
an element $\xi\in\G$ on some element $\eta\in {\rm End}\,{\cal H}$ 
is defined as follows \cite{BTR,Bie}:
\be{1.1}
  ({\rm ad}_q \xi) \;\eta \,=\, 
 \sum\limits_k \; {\xi}_k^1 \; \eta \;  S({\xi}_k^2), 
\ee
where $\xi^a_k$ are the components of the co-product
$\D\xi=\sum\limits_k {\xi}_k^1 \o {\xi}_k^2\in\G\o\G$ and 
$S(\xi)\in\G$ is the image of $\xi$ under action of the 
antipode.

{}From the physical point of view, the space $\H$ in (\ref{1.1}) is the
Hilbert space of the model in question. For $\G$ being a (quantum) Lie 
algebra it is often chosen as the corresponding  model space,
\hbox{$\M=\bigoplus_{I} {\cal H}_{I}$} ($I$ runs over all
highest weights and each subspace ${\cal H}_{I}$ appears with multiplicity 
one).

Let $\rho^J : \G\mapsto {\rm End}\ V^J$ be a highest weight $J$ irrep 
of $\G$ with the carrier vector space $V^J$ of the dimension $\d_J$.
The set of operators $\{T_m^J\}_{m=1}^{\d_J}$ acting on 
the Hilbert space $\cal H$ is called a {\em tensor operator} (of highest 
weight $J$) if
\be{1.2}
 ({\rm ad}_q \xi) \;T^J_m \,=\, 
 \sum\limits_n \,T^J_n \, \left(\rho^J(\xi)\right)_{nm}
 \;\;\;\; {\rm for \ all \ } \xi\in\G \ .
\ee
If $\{T_m^I\}_{m=1}^{\d_I}$ and 
$\{T_{m'}^J\}_{m'=1}^{\d_J}$ are tensor operators acting
on the same Hilbert space, then, using corresponding (deformed) 
Clebsch-Gordan coefficients, we can construct tensor operator of weight 
$K$ as follows \cite{BTR,Bie}:
\be{f}
 T^K_{m''} \,=\, \sum_{m,m'} \, \CG_q \, T^I_m \, T^J_{m'}  \ .
\ee
This formula describes {\it fusion} of tensor operators.

In the case of $\G=U_q(sl(2))$  tensor operators, 
$\{T_m^J\}_{m=-J}^{J}$, are labeled by spin $J$ and the 
definition (\ref{1.2}) acquires the form:
\be{1.3}
\ar{rcl}
 X^{\pm} \; T^{\,J}_m \, q^H - q^{H\mp 1} \; T^{\,J}_m \; 
 X^{\pm} &=& \sqrt{[J\mp m][J\pm m+1]}\, T^J_{m\pm 1} ,\\ [2mm]
 q^H \; T_m^{\,J} \; q^{-H} &=& 
 q^m \, T_m^{\,J} \ , \er
\ee
where $[x]=(q^{x}-q^{-x})/(q-q^{-1})$ is a $q$-number and
$X_\pm$, $H$ are the generators with the commutation relations
$$
 [ \,H\, , \, X_\pm\, ] \,=\, \pm X_\pm \ , \;\;\;\;
 [ \,X_+\, , \, X_-\, ] \,=\, [\,2\,H\,] \ .
$$
An example of (deformed) tensor operator (of spin 1) is provided by the 
following set of combinations of the generators:
\be{1.4}
 T^1_1 \,=\, q^{-H}\, X_+ \ , \;\;\; T^1_0 \,=\, 
 (q^{-1}\,X_-\,X_+ - q\,X_+\,X_-)/\sqrt{[2]} \ , \;\;\; 
 T^1_{-1} \,=\,-q^{-H}\, X_- \ . 
\ee 
Notice, however, that this is a rather special case because, generally 
speaking, components of tensor operators act on the model space as 
shifts between different subspaces $\H_I$, whereas for the components 
$T^1_m$ in (\ref{1.4}) $\H_I$ are invariant subspaces.

Let us remark that along with the tensor operator of {\em covariant} type
introduced in (\ref{1.2}) one can define a {\em contravariant} tensor 
operator as the set  of operators 
$\{{\ov{T}}_m^{\,J}\}_{m=1}^{\d_J}$ obeying the following relations:
\be{1.5}
({\rm ad}_q \xi)\;\ov{T}^{\,J}_m =\sum\limits_n \Bigl(\rho^J(S(\xi))
\Bigr)_{mn}\,\ov{T}_n^{\,J} \;\;\;\; {\rm for \ all \ } \xi\in\G \  .
\ee
{}Further we shall consider only the covariant case since the theory
and computations for the contravariant case are quite analogous. 

In the case of quasi-triangular Hopf algebra we can 
describe tensor operators using $R$-matrix language. 
Let $\rho^J$ be a highest weight $J$ irrep of $\G$ with 
the carrier space $V^J$ of dimension $\d_J$. Let 
$U^J\in{\rm End}\,V^J \o {\rm End}\,\H$ be a matrix obeying the following 
$R$-matrix relations:
\be{2.1}
 \up{1}{L}\^I_{\pm}\,\up{2}{U}\^{J} \,=\;
 \up{2}{U}\^{J} \, R^{IJ}_{\pm} \, \up{1}{L}\^I_{\pm} \ ,
\ee
where $L^I_\pm$ and $R^{IJ}_\pm$ are defined as in (\ref{0.1})-(\ref{0.2}).
Equations (\ref{2.1}) are equivalent \cite{Byt} to the statement
that each row of $U^J$ satisfies (\ref{1.2}); that is, all rows of $U^J$
are tensor operators of weight $J$.
We shall refer to $U^J$ as {\it generating matrices} because,
according to the Wigner-Eckart theorem (see, e.g., the comments 
in \cite{Byt}), matrix elements of entries of $U^J$ evaluated on 
vectors from $\H$ give the ($q$-deformed) Clebsch-Gordan coefficients.
Notice that if $U^J$ obeys (\ref{2.1}) and $M$ is a matrix with 
entries commuting with all the elements from $\G$, then
\be{TU}
  \wt{U}^J\,=\,M\,U^J
\ee
also obeys (\ref{2.1}), i.e., $\wt{U}^J$ also is a generating matrix.

The matrix $U^J$ in (\ref{2.1}) may have an arbitrary number of rows.
However, it is more natural to consider the case of $U^J$
being a square matrix; therefore, from now on we shall regard it as
$\d_J\times\d_J$ matrix.

Now let $U^I$ and $U^J$ be two generating matrices. The fusion formula 
(\ref{f}) can be written in the $R$-matrix language as follows \cite{Byt}: 
\be{2.2}
 U^{IJ}_{K} \,=\, P^{IJ}_{K} \, F^{IJ} \, 
 \up{2}{U}\^J \, \up{1}{U}\^I  \, P^{IJ}_{K} \ \  
 \in {\rm End}\,(V^I\o V^J) \o {\rm End}\,\H \ ,
\ee
where the l.h.s. is a new generating matrix of weight $K$ written in the 
basis of \hbox{ $V^I\o V^J$}. 
Here $F^{IJ}$ is an arbitrary $(\d_I\,\d_J)\times(\d_I\,\d_J)$ matrix 
whose entries commute with all elements of $\G$; and 
$P^{IJ}_{K}\in {\rm End}\,(V^I\o V^J)$ stands for the 
projector (i.e., $(P^{IJ}_{K})^2=P^{IJ}_{K}$) onto the subspace in 
$V^I\o V^J$ corresponding to the representation $\rho^K$ (cf. \S 3.2). 

One can rewrite (\ref{2.2}) in the conventional basis of the space $V^K$:
\be{2.3}
 U^K_{\,mn} \,=\, e_m^t\; U^{IJ}_K \; e_n , \;\;\;\;
 m,n \,=\, 1,.., \d_K. 
\ee
Here $\{e_n\}$ is an orthonormal set of the eigenvectors of the projector 
$P^{IJ}_K$; that is, $e^t_m\,e_n=\d_{mn}$ and
\hbox{$P^{IJ}_K=\sum_{n=1}^{\d_K} \! e_n\o e^t_n$}.

{}Formula (\ref{2.2}) resembles the fusion formula for $R$-matrices 
\cite{KS}:
\be{f1}
 \up{1,32}{R}\^{LK}_{\pm} \,=\, \up{23}{P}\^{IJ}_K\,
 \up{13}{R}\^{LJ}_{\pm}\,
 \up{12}{R}\^{LI}_{\pm}\, \up{23}{P}\^{IJ}_K \ ,
\ee
where the l.h.s. stands for $R^{\,LK}_{\,\pm}$ written in the basis of 
\hbox{$V^L\o V^I\o V^J$} and we use notations of \cite{KS}. Of course, the
origin of both eqs. (\ref{2.2}) and (\ref{f1}) is the Hopf structure of $\G$.

The fusion formulae given above are of direct practical use 
since they allow to construct corresponding objects (generating matrices
and $R$-matrices) for higher representations starting with those for 
the fundamental irreps. 

Let us now introduce the Clebsch-Gordan maps, 
${\rm C}[IJK]:V^I\o V^J\mapsto V^K$ and 
${\rm C}'[IJK]:V^K\mapsto V^I\o V^J$.
They are given by\footnote{
${\rm C}[IJK]$ and ${\rm C}'[IJK]$ in (\ref{2.6}) can be regarded as
rectangular matrices with numerical entries if vectors in $V^I$, 
$V^J$, $V^K$ are realized as usual numerical vectors. They act on these
vectors by matrix  multiplication from the left. For instance, 
the result of action ${\rm C}'[IJK]$ on a vector 
$a=\sum_n a_n e_n \in V^K$ is 
{${\rm C}'[IJK]\,a=\sum_n a_n \wh{e}_n\equiv \wh{a}$} --
the same vector but written in the basis of the tensor product 
$V^I\o V^J$. }
\be{2.6}
 {\rm C}[IJK] \,=\, \sum\nolimits_{n=1}^{\d_K} \! 
 \wh{e}_n\o e_n^t \ , \;\;\;\; 
 {\rm C}'[IJK] \,=\, \sum\nolimits_{n=1}^{\d_K} \! 
 e_n\o \wh{e}_n^{\,t}  \ ,
\ee 
where $\wh{e}_n$ stands for the vector $e_n$ rewritten in the basis 
of the space $V^K$ (cf. \S 3.2). The main properties of the CG maps are:
\be{2.7}
 {\rm C}[IJK]\,\left( (\rho^I\o \rho^J)\,\D(\,\xi\,) \right) \, 
 {\rm C}'[IJK]  \,=\, \rho^K(\xi) \ \ {\rm for\ any\ } \xi\in\G ,
\ee
\be{2.8}
 \sum\limits_K {\rm C}'[IJK] \,{\rm C}[LMK]
= \d_{IL}\,\d_{JM},\;\;\;\;\;
{\rm C}[IJK]\,{\rm C}'[IJL]=\d_{KL} \ .
\ee

With the help of the CG maps we can rewrite eqs. (\ref{2.3})-(\ref{f1}) in 
the following form:
\ba
 U^K &=& {\rm C}[IJK] \, F^{IJ} \, \up{2}{U}\^J \, \up{1}{U}\^I \, 
 {\rm C}\,'[IJK] \ ,  \label{2.4} \\ [1mm]
 R^{LK}_\pm &=& \up{23}{\rm C}[IJK] \, \up{13}{R}\^{LJ}_\pm\,
 \up{12}{R}\^{LI}_\pm\, \up{23}{\rm C}\!{}'[IJK] \ . \label{2.5}
\ea

\subsection{Exact generating matrices }
\def\x{\vec{\rm x}}

\red For generating matrices, as they have been defined above,
rows are not necessarily linearly independent tensor operators.
However, in the case of $\G=U_q(sl(n))$ (and very probably even in the
case of $U_q(\J)$ for any semi-simple $\J$) there exists a scheme which
allows to  obtain an example of generating matrix with all rows being
linearly independent tensor operators. Actually, this scheme has been 
developed in studies of the quantum WZNW model \cite{AF,CG,BS,BF}.
Let us now describe it (remember that we deal with the case of $|q|=1$).

Let $\J$ be a semi-simple Lie algebra of rank $n$. Introduce $n$-dimensional
vector $\p=2J+\rho$, where $J$ runs over all highest weights and $\rho$ 
is the sum of the positive roots of $\J$. Let $\C$ be a commutative algebra 
of functions on the weight space of $U_q(\J)$, i.e., an algebra of functions 
depending on the components of $\p$.

Next, let us introduce two auxiliary objects: 
\hbox{$D=q^{2\vec{H}\o\p}\in\G\o\C$} and
\hbox{$\Omega=q^{4\vec{H}\o\vec{H}}\in\G\o\G$}, where 
\hbox{$\vec{A}\o\vec{B}$} is understood as 
\hbox{$\sum_{i=1}^n A_i\o B_i$};
and $H_i$ are the basic generators of the Cartan subalgebra of $\G$. 
Finally, we define the homomorphism $\s:\C\mapsto \G\o\C$, such that
\be{3.1}
 \s(\p) \,=\, e\o\p + 2\vec{H}\o e .
\ee
Now we can look for objects $\R_\pm(\p)\in \G\o\G\o\C$ 
[obeying the standard relation $\R_-(\p)=(\R_+'(\p))^{-1}$, where $'$ means  
permutation of the first two tensor factors] 
which are solutions of the following equations:
\ba
 & \up{12}{\R}_\pm(\p) \, \up{13}{\R}_\pm(\p_2) \, 
 \up{23}{\R}_\pm(\p) \,=\,
 \up{23}{\R}_\pm(\p_1) \, \up{13}{\R}_\pm(\p) \, 
 \up{12}{\R}_\pm(\p_3) \ , & \label{3.2}  \\ [2mm]
 & [\, \R_\pm(\p) \, , \, q^{H_i}\o q^{H_i}\o e \, ] \,=\, 0 \ \ \ \ 
 {\rm for\ all}\ \  i \ ,& \label{3.3} \\ [1mm]
 & \R_-(\p) \, (e\o D) \;=\, (\Omega\o e) \, (e\o D) \, \R_+(\p) \ , &  
 \label{3.4} \\ [1mm] 
 & \R_\pm^{\,*} (\p) \,=\,  \R_\pm^{-1}(\p) \ . & \label{con}
\ea
The subscript $i=1,2,3$ of the argument of $\R_\pm(\p)$ means that
this argument is shifted according to (\ref{3.1}) and the $\vec{H}$-term 
appears in the $i^{th}$ tensor component.\footnote{
To make these shifts more transparent, let us introduce the element 
\hbox{$Q=e^{2\vec{H}\o\x}$}, where components of $\x$ are such that
\hbox{$[p_i,{\rm x}_j]=\d_{ij}$}. Then 
\hbox{$\up{12}{\R}_\pm(\p_3)=\up{3}{Q}\^{-1}\up{12}{\R}_\pm(\p)\up{3}{Q}$},
etc.\ Notice that the shifted matrices  belong to 
\hbox{$\G\o\G\o\G\o\C$}.\\ [-2.9mm] {} }
It is easy to verify that for $|q|=1$ eq. (\ref{con}) [unitarity of
$\R_\pm(\p)$ ] is consistent\footnote{
The conjugation of an object belonging to $m$-fold tensor product 
$\G^{\o m}$ is understood as follows: \\ [1mm]
\hbox{$(\xi_1\o\xi_2...\o\xi_m)^*=\xi_1^*\o\xi_2^*...\o\xi_m^*$}.}
with eqs. (\ref{3.2})-(\ref{3.4}).
Eq. (\ref{3.3}) is the same symmetry 
condition which is known\footnote{
Recall that the quantization $U(\J)\rar U_q(\J)$ does not deform
the co-multiplication for elements of the Cartan  subalgebra of $\J$.}
for the standard $R$-matrices of $U_q(\J)$.

In general, for $U_q(\J)$ there exists a family of solutions of 
eq. (\ref{3.2}). A solution $\R^{IJ}_\pm(\p)$ obeying, in addition,
the conditions (\ref{3.3})-(\ref{con}) is most remarkable among them
\cite{AF,CG,KR,BBB,BS}. Entries of such $\R^{IJ}_\pm(\p)$ coincide 
with the corresponding (deformed) $6j$-symbols.

Let us now consider an element $U\in \G \o {\rm End}\,\H$ which obeys the 
equations
\ba
 \R_\pm (\p) \up{2}{U} \,\up{1}{U} &=&  \up{1}{U} \, 
 \up{2}{U} R_\pm \ , \label{3.5} \\ [1mm]
 \up{1}{U} \, (e\o D) &=& 
 (\Omega\o e) \,  (e\o D) \, \up{1}{U}  \ , \label{3.6} \\ [1mm]
 U^{-1}\,D\,U &=& q^{C\o e}\,L_+\,L_-^{-1} \ , \label{L}
\ea
where $C\in\G$ and \hbox{$\rho^J(C)=2J(J+\rho)$}.
It can be shown \cite{AF,CG,BS,BF} that such $U$ 
(if it exists) is a generating matrix\footnote{ 
To be more precise, $U$ and $\R_\pm(\p)$ are not matrices but so called 
{\it universal} objects. \  If we fix  representations of 
their $\G$-parts: 
$U^J=(\rho^J\o id) U$, $\R_\pm^{IJ}(\p)=(\rho^I\o \rho^J )\R_\pm(\p)$,
we obtain a generating  matrix and  $\C$-valued counterparts of the
standard $R$-matrices.\\ [-2.9mm] {} }
for $\G=U_q(\J)$. Notice that (\ref{3.2})-(\ref{3.4}) are nothing but 
consistency conditions for (\ref{3.5})-(\ref{3.6}). 
The relation (\ref{3.6}) is a matrix form of the equation
\be{3.7}
        U \, \p \,=\, \s(\p) \, U .
\ee
Eq. (\ref{L}) plays the role of a normalization condition.

Observe that from the group of transformations (\ref{TU}) only the
following subgroup
\be{DU}
 U \mapsto  D^\a \, U \ , \ \  \a\in\Re \ 
\ee
survives for the solution $U$ of eqs. (\ref{3.5})-(\ref{L}).
Additionally, the rescaling  
\be{fU}
 U \mapsto (e\o f) \,U \, (e\o f)^{-1}
\ee
with an arbitrary element $f\in\C$ is allowed.

Let us explain why the generating matrix obeying (\ref{3.5})-(\ref{3.7}) 
is of particular interest from the point of view of the theory of tensor 
operators. Notice that the property (\ref{3.7}) ensures that rows of
$U^J=(\rho^J\o id)U$ are linearly independent tensor operators.\footnote{
To clarify this statement, we can rewrite eq. (\ref{3.7})
in the following form: $[U^J,p_i]=2H^J_i U^J$, $i=1,..,n$. In the 
conventional basis (where the generators $H_i$ of the Cartan subalgebra
are diagonal) the last relation for $k$-th row of $U^J$ becomes:
$[U^J_k,p_i]=2(H^J_i)_{kk} U^J_k$, $i=1,..,n$. Since the elements 
$H_i$ are linearly independent, we infer that the rows of $U^J$ generate 
linearly independent shifts on the space $\C$.} 
In other words, if for a given irrep $\rho^J$ and a given vector 
\hbox{$|I,m\rangle$} in the model space $\M$ of $\G$ we consider the 
set of vectors \hbox{$U^J_{ij}\,|I,m\rangle$} \ 
\hbox{$i,j=1,..,\d_J$}, then all non-vanishing vectors in 
this set are pairwise linearly independent. In particular, if $J$ is a 
fundamental representation, then entries of $U^J$ provide a set of basic 
shifts on $\M$. 
Thus, solutions of (\ref{3.5})-(\ref{L}) present very special 
but, in fact, the most interesting case of generating matrices. 
We shall call them {\it exact} generating matrices. 

An important from the practical point of view property of exact 
generating matrices is that the matrix elements  
\hbox{$\langle K,m''|U^J_{ij}|I,m'\rangle$} coincide 
up to some $p$-dependent factors allowed according to (\ref{DU})-(\ref{fU})
with the Clebsch-Gordan coefficients \hbox{$\CG_q$} 
(with the weights $I$, $J$, $K$ restricted by the triangle inequality).

Let us tell briefly about the physical content of the relations given 
above. Equation (\ref{3.2}) has appeared in various forms in studies of 
quantum versions of the Liouville \cite{CGN,Bab}, Toda \cite{Td} and 
Calogero-Moser \cite{ABB} models. 
In these models $\R(\p)$ is interpreted as a dynamical $R$-matrix.
{}From the point of view of the theory of tensor operators relations 
(\ref{3.1})-(\ref{3.7}) most closely connected with a quantization of the 
WZNW model \cite{AF,CG,BS}. Here $\R(\p)$ plays a role of the braiding
matrix and eqs. (\ref{3.5})-(\ref{L}) with appropriate dependence 
on the spatial coordinate (or its discretized version) describe vertex 
operators. Let us mention that in the WZNW theory
the quantum-group parameter of $\G=U_q(\J)$ is given by 
\hbox{$q=e^{i\g\hbar}$}, where \hbox{$\hbar\!>\! 0$} is the Planck 
constant and the deformation parameter \hbox{$\g>0$} is interpreted as a 
coupling constant. 
This provides the motivation to study the case of $|q|=1$.

\section{\nz\bf FUSION OF EXACT GENERATING MATRICES }
\setcounter{equation}{0}
\subsection{Formulation of the problem}

\red Suppose we are given two generating matrices, $U^I$ and $U^J$, for 
some irreps, $\rho^I$ and $\rho^I$, of $\G$. Then by formula (\ref{2.2}) 
we can build up a generating matrix $U^K$ for every irrep $\rho^K$ which 
appears in the decomposition of $\rho^J\o\rho^J$. For the sake of shortness 
we shall call them {\it descendant} matrices. However, as we have explained 
before, it is natural to deal not with all possible generating matrices 
but only with exact ones, i.e., with those which obey additional equations 
(\ref{3.5})-(\ref{3.7}) with $\R_\pm(\p)$, $D$ and $\Omega$ introduced 
above. Thus, if $U^I$ and $U^J$ are exact generating matrices, it is
natural to look for such a matrix $F^{IJ}$ that descendant generating 
matrices $U^K$ obtained by the formula (\ref{2.2}) would be also exact.

Let us underline that this problem would not arise if eq. (\ref{3.5}) 
contained in the l.h.s. the standard $R$-matrix instead of $\R(\p)$. 
Indeed, for an operator-valued matrix \hbox{$g^J\in{\rm End}\,V^J\o\G$} 
(it may be regarded as $L$-operator type object\footnote{
{}For $\G$ replaced by its dual $\G'$ the matrix $g^J$ is regarded as
the quantum group-like element. \ The  fusion formulae
(\ref{4.0b})-(\ref{4.0c}) are also valid in this case.\\ [-2.9mm] {} })
which obeys the usual quadratic relations\footnote{
We prefer this order of auxiliary spaces in (\ref{4.0a}) since it is the 
same as in (\ref{3.5}).}
\be{4.0a}
   R^{IJ} \, \up{2}{g}\^J \; \up{1}{g}\^I \,=\, 
   \up{1}{g}\^I \; \up{2}{g}\^J \, R^{IJ} \ ,
\ee
the fusion formula is well known [eq. (\ref{f1})
is its particular realization]:
\be{4.0b}
   (g^K)_{mn} \,=\, e_m^t\; \up{2}{g}\^J \, \up{1}{g}\^I \; e_n \ , 
\ee
where $e_n$, $n=1,..,\d_K$ are the eigenvectors of the 
projector $P^{IJ}_K$ introduced in \S 1.2. For example, in the case of 
$\G=U_q(sl(2))$, starting with $g^{\frac 12}$ and applying (\ref{4.0b}) 
iteratively, one obtains matrices $g^J$ for any spin $J$:
\be{4.0c}
\!\!\!\! g^0 = 1 \, , \;\;\; 
 g^{\frac 12}=\left(\ar{cc} a & b \\ c & d \er\right)  , \;\;\;
 g^1=\left(\ar{ccc} a^2 & q^{-\frac 12}\sqrt{[2]}\,ab & b^2 \\ 
 q^{-\frac 12}\sqrt{[2]} \, ac & ad + q^{-1}bc & q^{-\frac 12}\sqrt{[2]} \, 
 bd \\  c^2 & q^{-\frac 12}\sqrt{[2]} \, cd & d^2 \er\right) , \, \dots
\ee

{}For generating matrices the fusion problem is more complicated because 
$\R(\p)$ in eqs. (\ref{3.2}) and (\ref{3.5}) is an attribute not of Hopf 
algebra but of quasi-Hopf algebra. In this section we shall discuss some 
general aspects of the fusion problem in the quasi-Hopf case. In the next 
section we shall consider an example -- the case of $U_q(sl(2))$.

It should be also underlined that the fusion problem (as it is formulated 
above) does not appear if the language of {\it universal} objects 
(see, e.g., \cite{MS,BS}) is used instead of the language of 
operator-valued matrices. For example, instead of the set of matrices 
{$g^J\in{\rm End}\,V^J\o\G$} obeying (\ref{4.0a}) we could introduce 
the element $g\in\G\o\G$ and fix its {\it functoriality}
relation as follows:
\be{4.0d}
 (\D\o id) \, g \,=\, \up{2}{g} \, \up{1}{g} \ .
\ee
Then both quadratic relations (\ref{4.0a}) and fusion formula (\ref{4.0b})
can be obtained from (\ref{4.0d}) with the help of the axioms of 
quasi-triangular Hopf algebra. In fact, in this approach we actually do 
not need the fusion formula because each $g^J$ can be obtained simply by
evaluation of the universal element $g$ in the corresponding 
representation: $g^J=(\rho^J\o id)g$.

Similarly, we could introduce the universal object $U\in\G\o{\rm End}\,\H$
with the functoriality relation \cite{BS}:
\be{4.0e}
   (\D\o id) \, U \,=\, \F \, \up{2}{U} \, \up{1}{U} \ ,
\ee
where $\F$ obeys a certain set of axioms. Then quadratic relations 
(\ref{3.5}) [with $\R(\p)$ constructed from $\F$ and $R$ according to
(\ref{4.7})] would be consequences of (\ref{4.0e}).
Again, fixing representation of $\G$-part of the universal element $U$, 
we obtain a generating matrix $U^J=(\rho^J\o id)U$ and, therefore, we
do not need the fusion formula. 

Although the language of universal objects is more convenient in 
abstract theoretical constructions, in practice we usually do not have
explicit formulae for involved universal objects (or they are quite 
cumbrous; see, for instance, the universal $R$-matrices for $U_q(sl(n))$ 
\cite{uni}). 
Therefore, in the present paper we intentionally have adopted the matrix 
language to discuss how to construct exact generating matrices for an
arbitrary irrep from those for given irreps without invoking to
universal formulae.

\subsection{Quasi-Hopf features }

\red Let us remind that an associative algebra $\G$ equipped with 
co-product, co-unit and antipode is said to be 
quasi-Hopf algebra \cite{Dr2} if its co-multiplication is 
``quasi-coassociative''; that is, for all $\xi\in\G$ we have
\be{4.1}
  ((id\o\D) \, \D(\xi)) \, \Phi \,=\, \Phi \, ((\D\o id)\, \D(\xi)) \ ,   
     \ \ \  (\eps\o id)\,\D(\xi) \,=\, (id\o\eps)\,\D(\xi) \,=\, \xi.
\ee
Here $\Phi\in\G\o\G\o\G$ is an invertible element (the 
{\it co-associator}) which must satisfy certain equations. For 
quasi-triangular quasi-Hopf algebra it is additionally postulated
that there exists an invertible element $\R\in\G\o\G$ 
(the {\it twisted} $R$-matrix) such that
\ba
 \R\,\D(\xi)  \,=\, \D'(\xi)\,\R   && 
 {\rm for\ all}\ \xi\in\G \ , \label{4.2} \\ [1mm] 
 \!\!\!\! (\D\o id) \,\R = \Phi_{312} \up{13}{\R} {\Phi}^{-1}_{132}
 \up{23}{\R} \Phi_{123}  &\!\!\! ,& \!\!
 (id\o\D) \,\R = {\Phi}^{-1}_{231} \up{13}{\R} {\Phi}_{213}
 \up{12}{\R} {\Phi}^{-1}_{123} \ , \label{4.3} 
\ea
\be{4.4} 
 (\eps\o id)\, \R \,=\, (id\o\eps)\, \R \,=\, e \ .  
\ee
The analogue of Yang-Baxter equation for $\R$ follows from (\ref{4.2})
and (\ref{4.3}) and looks like 
\be{4.5}
 \up{12}{\R}\,{\Phi}_{312}\,\up{13}{\R}\,{\Phi}^{-1}_{132}\,
 \up{23}{\R}\,{\Phi}_{123} \,=\, {\Phi}_{321}\,\up{23}{\R}\, 
 {\Phi}^{-1}_{231}\,\up{13}{\R}\,{\Phi}_{213}\,\up{12}{\R} \ .
\ee

A crucial observation \cite{CG,BBB,BS} is that the construction 
of exact generating matrices, which we described in {\S} 1.3, involves the 
quasi-Hopf algebra $\G_\F$ (where $\R_+(\p)$ plays the role of the element 
$\R$) obtained as a twist of the symmetry algebra $\G$. 
More precisely, there exists an invertible element $\F(\p)\in\G\o\G\o\C$
such that one can construct with its help from standard co-multiplication
and $R$-matrices (which obey axioms of Hopf algebra) the following objects
which obey all the axioms of quasi-triangular quasi-Hopf algebra:\footnote{
In fact, here we deal with some generalization of the Drinfeld's
scheme, since $\F(\p)$, $\R_\pm(\p)$  and $\Phi(\p)$ possess 
additional $\C$-valued tensor component. But all Hopf-algebra operations 
are applied  only to $\G$-parts of these objects.}
\ba
 \D_\F(\xi) &=& \F^{-1}(\p) \,\D(\xi)\,\F(\p)\ \ \ {\rm for\ all}\ 
 \xi\in\G \ ,\label{4.6} \\ [1mm] 
 \R_\pm(\p)  &=& (\F'(\p))^{-1} \, R_\pm \,\F(\p) \ \ \in\G\o\G\o\C
 \ ,\label{4.7} \\ [1mm]
 \Phi(\p)_{123} &=& \up{12}{\F}\^{-1}(\p_3) \up{12}{\F}\!(\p) \ \
 \in\G\o\G\o\G\o\C \ . \label{4.8}
\ea
Here $\p_i$ denotes, as before, the shift (\ref{3.1}) with the 
$\vec{H}$-term appearing in $i^{th}$ tensor component. 
These formulae show that
equation (\ref{3.2}) introduced above is a particular realization of 
the abstract form (\ref{4.5}) of the twisted Yang-Baxter equation. 

The fact that $\R_\pm(\p)$ introduced in (\ref{3.2})-(\ref{3.4}) admit 
decomposition of type (\ref{4.7}) is crucial in the context of
the fusion problem for exact generating matrices. 
Indeed, suppose we are given two exact 
generating matrices, $U^I$ and $U^J$, which obey (\ref{3.5})-(\ref{3.7}) 
with certain $\R^{IJ}_\pm(\p)$. Applying the formula (\ref{2.2}) with 
some matrix $F^{IJ}(\p)$ (it may be $\C$-valued) to these $U^I$ and 
$U^J$, we get new matrix $U^{IJ}_K$. It automatically obeys (\ref{2.1}). 
Moreover, it is easy to verify that exchange relations between 
$U^{IJ}_K$ and any exact generating matrix $U^L$ have the form (\ref{3.5}) 
but contain matrices $\up{1,32}{R}\^{LK}_\pm$ given by (\ref{f1})
in the r.h.s. and some $p$-dependent $R$-matrices in the l.h.s. The latter
look like following 
\be{4.9}
  \up{1,32}{\R}\^{LK}_\pm(\p) \,=\; 
  \up{23}{P}\^{IJ}_K \, \up{23}{F}\^{IJ}(\p_1) \,
  \up{13}{\R}\^{LJ}_\pm(\p) \, \up{12}{\R}\^{LI}_\pm(\p_3) \,
  (\up{23}{F}\^{IJ}(\p))^{-1} \, \up{23}{P}\^{IJ}_K \  ,
\ee
where, similarly as in (\ref{f1}), the basis of $V^L\o V^I\o V^J$
is used.
This is an analogue of the fusion formula (\ref{f1}) for standard 
$R$-matrices. 
The demand that the new generating matrix $U^{IJ}_K$ is exact, i.e.,
in particular, it obeys (\ref{3.5}), implies that expression (\ref{4.9}) 
rewritten in the basis of {$V^I\o V^J$} must coincide with 
$\R^{LK}_\pm(\p)$. Taking into account that $\R_\pm(\p)$ satisfies 
(\ref{4.7}) with some element $\F\in\G\o\G\o\C$, we get the equality
\ba
 && \R^{LK}_\pm(\p) \,=\, \up{23}{\rm C}[IJK] \,\up{23}{F}\^{IJ}(\p_1)\,
 \up{13}{\R}\^{LJ}_\pm(\p) \, (\up{21}{\F}\^{IL}(\p_3))^{-1} \, 
 \up{21}{\F}\^{IL}(\p)\, \times \label{4.9b} \\ [1mm] 
 && \times \,
 \up{12}{\R}\^{LI}_\pm(\p)\, (\up{12}{\F}\^{LI}(\p))^{-1} \, 
 \up{12}{\F}\^{LI}(\p_3)\,
 (\up{23}{F}\^{IJ}(\p))^{-1} \, \up{23}{\rm C}\!{}'[IJK] \ . \nn
\ea
The latter is equivalent due to (\ref{2.7}) and (\ref{4.6}) to the 
identity
$$ \ar{c}
 \left(\rho^L\o\rho^I\o\rho^J\right) \, (id\o\D_\F\!)\,\R \,=\, 
 (\up{23}{F}\^{IJ}(\p))^{-1} \,
 \up{23}{F}\^{IJ}(\p_1)\, \up{13}{\R}\^{LJ}_\pm(\p)\, \times \\ [1mm]
 \times \, (\up{21}{\F}\^{IL}(\p_3))^{-1} \, \up{21}{\F}\^{IL}(\p)\, 
 \up{12}{\R}\^{LI}_\pm(\p)\, 
 (\up{12}{\F}\^{LI}(\p))^{-1} \, \up{12}{\F}\^{LI}(\p_3) \ , \er  
$$
which, as we see from (\ref{4.3}) and (\ref{4.8}), takes place only if 
$F^{IJ}(\p)=\F^{IJ}(\p)$. 

\subsection{Properties of the twisting element  }

\red Let us give a r\'esum\'e of the previous paragraph. Let $U$ be
a universal generating matrix for a given symmetry algebra $\G$ 
(in the sense of \S 1.3) and let $\R_\pm(p)$ be the corresponding (twisted) 
$R$-matrix. Let $\F(\p)$ be the twisting element which transforms the Hopf 
algebra $\G$  into the quasi-Hopf algebra, $\G_\F$, for which $\R_\pm(\p)$ 
is an $R$-matrix (in the sense of eqs. (\ref{4.1})-(\ref{4.5})).  
If we are given two concrete representations of the exact generating matrix, 
$U^I$ and $U^J$ [and, hence, we know $\R^{IJ}_\pm(\p)$], then to construct 
its another representation $U^K$ we must substitute the matrix 
$\F^{IJ}(\p)=(\rho^I\o\rho^J)\F(\p)$ into the fusion formula (\ref{2.2}). 
An obstacle to application of this prescription is that explicit universal 
expressions for $\F(\p)$ and $\R(\p)$ are usually unknown. However, one
can formulate some conditions which $\F^{IJ}(\p)$ has to obey:\\

{\it 1.}\ {$\F^{IJ}(\p)$ is a solution of the following equation
 (for given $\R^{IJ}_\pm(\p)$):
\be{1}
   R^{IJ}_\pm \,\F^{IJ}(\p) \,=\, (\F^{IJ}(\p))' \, \R^{IJ}_\pm(\p) \ ;
\ee

{\it 2.}\ $\F^{IJ}(\p)$ obeys the symmetry condition
\be{2} 
 [\,\F^{IJ}(\p) , \, q^{H^I_i}\o q^{H^J_i}\o e \, ] \,=\, 0 \ \ \ \  
 {\rm for\ }\ i=1,..,n \ ;
\ee}

{\it 3.}\ $\F^{IJ}(\p)$ is such that all the entries of the matrix
\be{3a}
 U^{IJ}_0 \,=\,  P^{IJ}_{0} \, \F^{IJ}(\p) \,\up{2}{U}\^J \up{1}{U}\^I \, 
 P^{IJ}_{0} \ \ \in {\rm End}(V^I\o V^J)\o {\rm End}\H \ , 
\ee
commute with all entries of generating matrix $U^M$ for any weight $M$; \\ 

{\it 4.}\ {$\chi^{IJ}=\F^{IJ}(\p)(\F^{IJ}(\p))^*$ \ 
 is a $p$-independent object. \\ [1mm]

Let us comment these conditions. First of them is equivalent to eq.
(\ref{4.7}); its necessity has been explained in the previous paragraph. 
However, this is not a sufficient condition since, in general, (\ref{1}) 
possesses a family of solutions. 
In principle, we could select the right solution in
this family verifying whether a substitution of this solution in fusion 
formulae (\ref{4.9}) or (\ref{4.9b})  yields  matrices $\R(\p)$ obeying
(\ref{3.2})-(\ref{con}). But such a verification would be quite cumbrous
in practice. 

The second condition ensures that the descendant matrix $U^K$ obeys eq. 
(\ref{3.7}) and, as a consequence, eq. (\ref{3.6}). This can be easily 
checked applying (\ref{3.7}) to (\ref{2.2}). In fact, eq. (\ref{2}) 
implies that for our specific example of quasi-Hopf algebra 
the co-multiplication on the Cartan subalgebra is not deformed, i.e., 
it is the same as for $U_q(\J)$ and $U(\J)$.

The third condition is derived from eq. (\ref{4.4}) and the property 
$(\eps\o id)R=(id\o\eps)R=e$\/ of the standard $R$-matrices
(recall that $\eps$ stands for the trivial one-dimensional representation 
of $\G$). Indeed, applying $(\eps\o id)$ or $(id\o\eps)$ to (\ref{3.5}), 
we conclude that {$U^0=(\eps\o id)U \in {\rm End}\H$ } commutes with 
all entries of $U^J$ for any $J$. For $U^{IJ}_0$ in eq. (\ref{3a}) we
have: $U^{IJ}_0=U_0 P^{IJ}_0$.\footnote{
{} Hence $tr(U^{IJ}_0)$ coincides (up to a constant) with $U^0$.
The latter can be regarded \cite{Byt} as a generalization of the quantum 
determinant \cite{FRT,KS}. Notice also that
it can be written with the help  of the CG-maps as follows:
$ U^0 =  {\rm C}[IJ0] \, \F^{IJ}(\p) \, 
 \up{2}{U}\^J \up{1}{U}\^I \, {\rm C}'[IJ0]$ 
(assuming that $\rho^0\equiv\eps$ appears in the 
decomposition of $\rho^I\o\rho^J$).
See \S 3.2 for the case of $U_q(sl(2))$.}
Therefore, if the trivial representation
$\rho^0\equiv\eps$ appears in the decomposition of the product
$\rho^I\o\rho^J$, then the condition {\it 3} is non-trivial.

The fourth condition claims, in fact, that the element 
$\chi=\F(\p)(\F(\p))^*$ belongs to $\G\o\G$ (or, more precisely, that 
the last tensor component of $\chi$ in $\G\o\G\o\C$ is trivial).
To clarify this, let us first recall that for $|q|=1$
the standard co-multiplication has the following property with respect 
to the conjugation in $\G$: {$\D^*(\xi)=\D'(\xi^*)$}. On the other hand,
relations (\ref{con}) imply that (see also \cite{BS})
\be{*D}
 (\D_\F(\xi))^* \,=\, \D_\F(\xi^*) \ \ {\rm for\ all}\ \xi\in\G\ .
\ee
Next, we observe that eqs. (\ref{con}), (\ref{4.3}) and
(\ref{*D}) ensure a unitarity of the co-associator:
\be{*P}
        \Phi^*(\p)_{123} \,=\, \Phi^{-1}(\p)_{123} \ .
\ee
According  to (\ref{4.8}), the latter equation leads to the condition
{$\up{12}{\chi}(\p)=\up{12}{\chi}(\p_3)$} or, equivalently, to
$\up{3}{U}(\up{12}{\chi}(\p))(\up{3}{U})^{-1}=\up{12}{\chi}(\p)$.
This implies $p$-independence of $\chi$.

Actually, the element $\chi$ plays an important role in the theory of
exact generating matrices.  Let us mention here that due to eq. (\ref{con}) 
it satisfies the following relation:
\be{ch}
R_\pm \, \chi \,=\, \chi' \, R_\mp^{-1} \ .
\ee
We shall discuss some other properties of $\chi$ below, in \S 3.4.
  
\section{\nz\bf $\bf U_q(sl(2))$ CASE}
\setcounter{equation}{0}

\red In this section the preceding discussion will be illustrated 
on some explicit calculations. Although solutions for the twisted 
Yang-Baxter equation (\ref{3.2}) are known \cite{CG,Is,CGN,Td}
for the fundamental representations of $\G=U_q(sl(n))$, 
we shall consider here only the case of $U_q(sl(2))$.
But let us stress that in the more general case of $U_q(sl(n))$ the
computations would be essentially the same.

\subsection{$\R(p)$ and $U$ in the fundamental representation  }

\red As we mentioned above, the twisted Yang-Baxter equation possesses a 
family of solutions even in the simplest case of the fundamental 
representation of $U_q(sl(2))$. But imposing additional conditions 
(\ref{3.3})-(\ref{con}), we get the unique solution $\R_\pm(p)$ 
\cite{AF,CG,Is} which depends on the single variable $p=2J+1$ with $J$ 
being the spin. In fact, $\R_-(p)$ coincides with $\R_+(p)$ if $q$ is
replaced with $q^{-1}$.
The entries of $\R_\pm^{IK}(p)$  coincide with values of the $6j$-symbols 
involving spins $I$, $K$ and $J$ (see \cite{KR} for details). 
Moreover, the asymptotics of $\R_\pm(p)$ in the formal limits 
$q^p\rar+\infty$ and $q^{-p}\rar+\infty$ are given by:
\be{asR}
 \ar{ccl}
 \R_{\,\pm}(p) & \rar  & R_\pm \ \ \ \
 {\rm when}\ q^p\rar+\infty \ ,\\ [1mm]
 \R_{\,\pm}(p) & \rar & R_\mp^{-1} \ \ \
 {\rm when}\ q^{-p}\rar+\infty \ . \er 
\ee
That is, in these limits we return to the case of Hopf algebra;
in particular, the co-associator becomes trivial:
$\Phi_{123} \rar e\o e\o e$.
Furthermore, relations (\ref{asR}) together with eqs. (\ref{4.7}), 
(\ref{ch}) allow to add the following condition to the properties of 
$\F(\p)$ listed in \S 2.3:\\

{\it 5}.\ (Asymptotic behaviour)
\be{asF} 
 \F(q^p\rar+\infty) \,=\, e\o e \ , \;\;\;\;
 \F(q^{-p}\rar+\infty) \,=\, \chi \ , 
\ee
where the element $\chi\in\G\o\G$ was described in \S 2.3. 
Let us stress that we have derived this additional condition, in its 
present form, only for $\J=sl(2)$. It would be interesting to find  
analogues of (\ref{asF}) in the case when $\p$ has several components. \\

In the simplest non-trivial case, $I=K=\frac 12$, the solution 
$\R^{IK}_\pm(p)$ is given by (all non-specified entries are zeros)
\be{5.1}
 \R_{\,+}^{\frac 12 \frac 12}(p) \,=\, 
 (\R_{\,-}^{\frac 12 \frac 12})' \,=\,
 q^{-1/2} \left( \ar{cccc} q&&& \\
 &  \frac{\sqrt{ [p+1] [p -1] }}{[p]} &  
 \frac{q^{p}}{[p]} &  \\ [2mm]
 &  -\frac{q^{-p}}{[p]} & 
 \frac{\sqrt{ [p +1] [p -1] }}{[p]} &  \\
 &&& q  \er \right) \ . 
\ee
Here $[x]$ stands, as usually, for $q$-number. 

Now let us turn to the solution $U^{\frac 12}$ of eqs. 
(\ref{3.5})-(\ref{L})
for $\R_{\,\pm}^{\frac 12 \frac 12}(p)$ given by (\ref{5.1}).
It was considered in different contexts in \cite{AF,CG,BF}
and has been shown to be unique up to the transformations 
(\ref{DU})-(\ref{fU}). In agreement with the general description 
in \S 1.3, the entries of 
\hbox{$U^{\frac 12}=\left(\ar{cc}U_1&U_2\\U_3&U_4\er\right)$} 
act on the model space 
\hbox{$\M=\oplus_{J=0}^\infty\H_J$} as the basic shifts (see Fig.1).

A particular realization of $U^{\frac 12}$ can be written (see, e.g., 
\cite{BF}) in terms of multiplication and difference derivative
operators for two complex variables:
\be{5.2}
 U^{\frac 12} \,=\,
 \left(\ar{cccc}  z_1\,q^{\frac 12 z_2\partial_2} & 
 z_2\,q^{-\frac 12 z_1\partial_1} \\ 
 -z_2^{-1}[z_2\partial_2]\,q^{-\frac 12 (z_1\partial_1+1)} &
 z_1^{-1}[z_1\partial_1]\,q^{\frac 12 (z_2\partial_2+1)}  \er\right)\,
 \frac{1}{\sqrt{[p]}} \ ,
\ee
where $p=z_1\partial_1+z_2\partial_2+1$. Entries of (\ref{5.2})
are operators of basic shifts on the model space realized as 
$D_q(z_1,z_2)$ -- the space of holomorphic functions of two 
complex variables equipped with such a scalar product (a deformation of 
the standard one) that the monomials  
\hbox{$|J,m\rangle:= \frac{z_1^{J+m}z_2^{J-m}}{\sqrt{[J+m]![J-m]!}}$} 
form an orthonormal basis. 

Moreover, the specific realization (\ref{5.2}) of the exact generating 
matrix of spin 1/2 may be called {\it precise} in the following sense. 
Matrix elements of its entries \hbox{$\langle J',m'|U_i|J'',m''\rangle$} 
evaluated on $D_q(z_1,z_2)$ exactly coincide with the 
CG coefficients 
\hbox{$\Bigl\{\ar{c}\!\st J'\;\:\;\;1/2\;\;\;J'' \\ [-4pt] 
\st m'\;\pm 1/2\;\:\,m'' \er\!\Bigr\}_q$} (four of them are 
non-vanishing); see \cite{BF} for more detailed comments. 
\begin{figure}
\begin{center}
\begin{picture}(200,150)(0,0)
\put(20,20){\vector(1,0){150}} \put(70,20){\vector(0,1){131}}

\multiput(90,20)(0,10){12}{\line(0,1){5}}
\multiput(110,20)(0,10){12}{\line(0,1){5}}
\multiput(130,20)(0,10){12}{\line(0,1){5}}

\multiput(70,80)(10,0){8}{\line(1,0){5}}
\multiput(70,100)(10,0){8}{\line(1,0){5}}
\multiput(70,120)(10,0){8}{\line(1,0){5}}

\put(73,147){$J$} \put(173,23){$m$}

\put(72,5){$m\!-\!\frac{1}{2}$} \put(106,5){$m$}
\put(122,5){$m\!+\!\frac{1}{2}$}

\put(42,78){$J\!-\!\frac{1}{2}$} \put(42,98){$J$}
\put(42,118){$J\!+\!\frac{1}{2}$}

\put(75,65){$U_4$} \put(75,125){$U_2$} 
\put(135,65){$U_3$} \put(135,125){$U_1$}

\thicklines
\put(110,100){\vector(1,1){20}} \put(110,100){\vector(-1,1){20}}
\put(110,100){\vector(1,-1){20}} \put(110,100){\vector(-1,-1){20}}

\put(-35,-25){Fig.1 \ Action of the operators $U_i$ on the model space.}
\end{picture}
\end{center} 
\end{figure}

Now we encounter the simplest version of the fusion problem --
to build up the exact generating matrix $U^1$ (of spin 1) 
from $U^{\frac 12}$. For this purpose we have to find an explicit form of 
the corresponding twisting element $\F(\p)$ in the fundamental 
representation.

Before going into the computations let us mention that a universal 
formula (i.e., applicable for representations of any spin) for solution 
$\wt{\R}(p)$ of eq. (\ref{3.2}) and a universal expression for 
$\wt{\F}(p)$ obeying (\ref{4.7}) with this $\wt{\R}(p)$ have been obtained 
in \cite{Bab}. But this solution $\wt{\R}(p)$ does not satisfy 
(\ref{3.4})-(\ref{con}) and, therefore, being evaluated, say, in the 
fundamental representation it differs from (\ref{5.1}). 
Thus, solution $\wt{U}$ of (\ref{3.5}) for such 
$\wt{\R}(p)$ would not be an exact generating matrices in our sense. 
In particular, the solution for spin 1/2 would differ from the
one given by (\ref{5.2}) and, therefore, would not have the remarkable 
properties mentioned above. 

Although, let us stress that such $\wt{U}$ still would be a generating 
matrix in the sense of definition (\ref{2.1}). Therefore, one could examine 
whether it can be converted into an exact generating matrix by means of 
the transformation $\wt{U}=M(p)U$ with \hbox{$M(p)\in\G\o\C$}.
If such $M(p)$ exists, then the following relations hold:
\ba
 \wt{\F}(p) &=& \Bigl( (\D\o id)M(p)\Bigr) \;\F(p)\,
  \Bigl( \up{1}{M}\!(p_2)\up{2}{M}\!(p) \Bigr)^{-1} 
 \ , \nn \\  \wt{\R}_\pm(p) &=& \up{2}{M}\!(p_1)\up{1}{M}\!(p)\;
 \R_\pm(p)\,(\up{1}{M}\!(p_2)\up{2}{M}\!(p))^{-1} \ , \nn
\ea
and we can construct our $\F(p)$ from $\wt{\F}(p)$ and $M(p)$. However, 
bearing in mind possible applications in the cases where no universal 
formulae for $\R(\p)$ are known, it is more instructive to give
a direct computation of $\F(p)$.

\subsection{Computation of $\F^{\frac 12 \frac 12}(p)$ and $U^1$}

\red The matrix $\F^{\frac 12 \frac 12}(p)$ must satisfy the 
conditions listed in \S 2.3. First of all, it must be a solution of 
the equation (\ref{1}), where $\R_\pm(p)$ in the r.h.s. are given by 
(\ref{5.1}) and the standard $R$-matrices in the l.h.s. are 
\be{6.0}
 R_{\,+}^{\frac 12 \frac 12} \,=\, q^{-1/2}\,\left(\ar{cccc} 
 q & & & \\ & 1 & \om & \\ &  & 1 & \\ & & & q \er\right) \ , \;\;\;
 R_{\,-}^{\frac 12 \frac 12} \,=\, q^{1/2}\, \left(\ar{cccc} 
 q^{-1} & & & \\ & 1 &  & \\ & -\om  & 1 & \\ & & & q^{-1} \er\right) \ ,
\ee
with $\om=q-q^{-1}$.
The symmetry condition {\it 2}\ dictates to look for the solution of
eq. (\ref{1}) in the following form:
\be{6.1}
 \F(p)= \left(\ar{cccc} 1 & & & \\ & \a (p) & \b (p) & \\ & \g (p) &
 \d (p) & \\ & & & 1 \er\right) \ .
\ee
The straightforward check shows that only two of the functions 
$\a(p)$, $\b(p)$, $\g(p)$, $\d(p)$ are independent, and we can express, 
say, the entries in the third row in (\ref{6.1}) via the entries of the 
second row. The result reads as follows:   
\be{6.2}
 \ar{c} \g(p)\,=\,-\frac{q^{-p}}{[p]}\,\a(p)+
 \frac{\sqrt{[p+1][p-1]}}{[p]}\,\b(p)\ , \;\;\;
 \d(p)\,=\,\frac{\sqrt{[p+1][p-1]}}{[p]}\,\a(p)+
 \frac{q^p}{[p]}\,\b(p)\ . \er
\ee

Now we shall use the condition {\it 3}. For this purpose, we can employ 
the the following formulae for the fundamental $R$-matrices of 
$U_q(sl(n))$ (see \cite{FRT} for details): 
\be{6.4}
 P_\pm \,=\, \ds \frac{q^{ \frac 1n \pm 1}\,\wh{R}_+ - 
 q^{-\frac 1n \mp 1}\,\wh{R}_- }{q^2-q^{-2}} \ ,
\ee
where $\wh{R}_\pm=PR_\pm$ [$P$ is the permutation matrix, i.e., {$PaP=a'$} 
for $a\in\G\o\G$] and $P_+$, $P_-$ are the 
projectors in ${\Ce}^n\o {\Ce}^n$ ($q$-symmetrizer and 
$q$-antisymmetrizer) of ranks $\frac{n(n+1)}{2}$ and $\frac{n(n-1)}{2}$, 
respectively. In the case of $U_q(sl(2))$ these projectors are given by
\be{6.5}
  P_+ \,=\, \left(\ar{cccc} 1&&&\\ &q^{-1}\,\la & \la& \\
  & \la &q\,\la &\\ &&&1 \er\right) \ , \;\;\;\;\;
  P_- \,=\, \left(\ar{cccc} 0&&& \\ &q\,\la & -\la& \\
  &-\la & q^{-1}\,\la &\\ &&&0 \er\right) \ , 
\ee
where $\la=\frac{1}{[2]}=(q+q^{-1})^{-1}$. It is easy to find their 
eigenvectors $\vec{x}_i$ such that $\vec{x}^{\,t}_i\vec{x}_j=\!\d_{ij}$:
\be{6.6}
P_+=\sum_{i=1}^3 \vec{x}_i\o\vec{x}^{\,t}_i \ , \;\;\;
\vec{x}_1 = \left(\ar{c} 1\\ 0\\ 0\\ 0\er\right),\;\;\;
\vec{x}_2=\sqrt{\la}\left(\ar{c} 0\\ q^{-1/2}\\ q^{1/2}\\ 0\er\right),\;\;\;
\vec{x}_3 = \left(\ar{c} 0\\ 0\\ 0\\ 1\er\right);
\ee
\be{6.7}
P_-=\vec{x}_0\o\vec{x}^{\,t}_0 \ , \;\;\;\;
\vec{x}_0=\sqrt{\la}\left(\!\ar{c} 0\\ q^{1/2}\\ -q^{-1/2}\\ 0\er\!\right).
\ee
According to (\ref{2.6}) we can construct from these vectors the following 
CG maps
\be{6.8} 
 \ar{c}
 \!\!\! {\rm C}[\frac 12 \frac 12 0] \,=\, \sqrt{\la}\,
 (0,\,q^{1/2},\, -q^{-1/2},\,0) \ ,\;\;\;
 {\rm C}[\frac 12 \frac 12 1] \,=\, \left(\ar{cccc} 1&0&0&0\\ 
 0&\sqrt{\la}\,q^{-1/2} & \sqrt{\la}\,q^{1/2}&0\\ 0&0&0&1 \er\right) \ .\er 
\ee
Now, substituting (\ref{6.5}) into (\ref{3a}), we can 
compute $U^0$. To be able to use the condition {\it 3} we have to compare 
$U^0$ with the central element of the algebra $\cal U$ generated by the 
entries $U_i$ of the matrix $U^{\frac 12}$ and the spin operator $p$. 
It has been shown in \cite{BF} that the only
nontrivial central element of the algebra $\cal U$ is given by the
following ``$p$-deformed'' determinant of $U^{\frac 12}$:
\be{6.9}
 \ar{c} {\rm Det}\, U^{\frac 12} \,=\, \left(U_1 U_4 - q\,U_2 U_3\right)
 \sqrt{\frac{[p]}{[p-1]}} \,=\,\left(q\,U_4 U_1 - U_3 U_2\right) 
 \sqrt{\frac{[p]}{[p+1]}} \ . \er
\ee
Omitting simple calculations, we give the result: $U^0$ coincides (up to
a numerical factor) with (\ref{6.9}) only if the constraint
\be{6.10}
 \a(p)\,\sqrt{[p+1]} - \b(p)\,\sqrt{[p-1]}\,=\,\ve\,\sqrt{[p]}.
\ee
holds. Thus, (\ref{6.1}) contains only one independent function.
Moreover, from the condition {\it 5} we infer that the
numerical constant $\ve$ in the r.h.s. of (\ref{6.10}) is fixed:
$\ve=q^{1/2}$.

{}Finally, we can use the conditions {\it 4} and {\it 5}. To apply the
former in practice, we can first consider the non-deformed case ($q=1$) 
when entries of $\F(p)$ are self-conjugated and then extend the solution
to generic $q$ in such a way that the condition {\it 5} would be satisfied. 
After simple calculations we get
\be{6.11}
 \ar{c}
 \a(p) \,=\, \d(p) \,=\, \frac{1}{[2]} \Bigl(\,q^{\frac 12}\,
 \sqrt{\frac{[p+1]}{[p]}} + q^{-\frac 12} \, \sqrt{\frac{[p-1]}{[p]}}
 \,\Bigr) \ , \\ [3mm]
 \b(p) \,=\, -\g(p) \,=\, {\frac{1}{[2]}} \Bigl(\,q^{-\frac 12}\,
 \sqrt{\frac{[p+1]}{[p]}} - q^{\frac 12}\, \sqrt{\frac{[p-1]}{[p]}}
 \,\Bigr)  \ . \er
\ee
Thus $\F^{\frac 12 \frac 12}(p)$ is found. Observe that \  
${\rm det}\,\F^{\frac 12 \frac 12}(p)=1$.

Now, substituting (\ref{6.8}) into (\ref{2.4})
and exploiting the explicit form of $\F^{\frac 12 \frac 12}(p)$, 
we can build up the exact generating matrix of spin 1. 
It looks like following:
\be{6.13}
 U^1 \,=\, \left(\ar{ccc} U_1^2 & q^{-\frac 12}\sqrt{[2]}\,U_1 U_2 
 & U_2^2 \\ [1mm] \sqrt{\frac{[2][p]}{[p+1]}} \, U_1 U_3 & 
 \sqrt{\frac{[p]}{[p+1]}} \, (q^{\frac 12}\, U_1 U_4 + q^{-\frac 12}\, 
 U_2 U_3) & \sqrt{\frac{[2][p]}{[p+1]}} \, U_2 U_4 \\ [1mm]
 U_3^2 & q^{-\frac 12}\sqrt{[2]}\,U_3 U_4 & U_4^2 \er\right) \ .
\ee

Let us briefly comment this formula. First, as we could expect [due to
eqs. (\ref{asF})], in the formal limit \hbox{$q^p\rar+\infty$} the 
structure of $U^1$ becomes identical to that of $g^1$ given in (\ref{4.0c}). 
Also, it is easy to see that the second row of (\ref{6.13}) 
can be identified (up to rescaling by $\sqrt{\frac{[2]}{[p+1][p-1]}}$)
with the spin 1 tensor operator (\ref{1.4}) constructed from the generators 
of $U_q(sl(2))$. Finally, the entries $U^1_{ij}$, $i,j=1,2,3$ act on the  
model space $\M$ as shifts from the state $|J,m\rangle$ to the states
$|J+(2-i),m+(2-j)\rangle$, which is natural because we have applied the 
fusion scheme to the matrix $U^{\frac 12}$ whose entries are basic shifts 
on $\M$.

{}Furthermore, if we substitute in (\ref{6.13}) the realization (\ref{5.2})
of the operators $U_i$, then $U^1$ also will be precise in the described
sense. Namely, it can be checked then that matrix elements
\hbox{$\langle J',m'|U^1_{ij}|J'',m''\rangle$} evaluated on $D_q(z_1,z_2)$
precisely coincide with the CG coefficients 
\hbox{$\Bigl\{\ar{c}\st \!J'\;\;\;\;1\;\;\;\;J''\\ [-4pt] 
\st m'\;\; 2-j \;\,m'' \er\!\Bigr\}_q$} (nine of them are non-vanishing). 
Thus, the fusion procedure preserves the ``preciseness'' of exact 
generating matrices, which is useful in practical applications.

\subsection{Another construction for $\F$}

\red The computations performed in the previous paragraph inspire us 
to introduce $p$-dependent counterparts of the projectors $P_\pm$ used
above. Indeed, we can consider the following analogue of
the decomposition formula (\ref{6.4}):
\be{7.1}
 \P_\pm \,=\, \ds \frac{q^{ \frac 1n \pm 1}\,\wh{\R}_+(\p) - 
 q^{-\frac 1n \mp 1}\,\wh{\R}_-(\p) }{q^2-q^{-2}} \ .
\ee
It is obvious from the formula 
\hbox{$\wh{\R}_\pm(\p)=(\F(\p))^{-1}\wh{R}_\pm\F(\p)$} that the 
objects $\P_+$ and $\P_-$ are also projectors of ranks $\frac{n(n+1)}{2}$ 
and $\frac{n(n-1)}{2}$, respectively. \footnote{
Notice that $P_-\F(p)\P_+=P_+\F(p)\P_-=0$. This is an alternative
form of eq. (\ref{1}) [for the case of fundamental representation].}
In the case of $U_q(sl(2))$ we find (cf. formulae (\ref{6.5}))
\be{7.2}
 \P_\pm \,=\, \left(\ar{cccc} \frac{1\pm 1}{2}&&&\\ [-0.5mm] 
 &\la\,\frac{[p\mp 1]}{[p]}& \pm\la\,\frac{\sqrt{[p+1][p-1]}}{[p]}&\\ [2mm]
 &\pm\la\,\frac{\sqrt{[p+1][p-1]}}{[p]} &
 \la\,\frac{[p\pm 1]}{[p]}&\\ &&& \frac{1\pm 1}{2} \er\right) \ . 
\ee
\def\x{\vec{\rm x}}

Repeating the procedure described in the previous paragraph, we can
find the eigenvectors $\x_i$ such that \hbox{$\x_i^{\,t}\x_j=\d_{ij}$}
and $\P_-=\x_0\o\x_0^{\,t}$,  $\P_+=\sum_{i=1}^3\x_i\o\x_i^{\,t}$.
Next, using the same formulae (\ref{2.6}), we can construct $p$-dependent
counterparts of the CG-maps. They look like following:
\be{7.3}
 \ar{ccl}
  {\rm C}_p[\frac 12 \frac 12 0] &=& \sqrt{\la}\,
 (0,\,\sqrt{\frac{[p+1]}{[p]}},\, -\sqrt{\frac{[p-1]}{[p]}},\,0) \ ,\\ 
 [2.5mm]      {\rm C}_p[\frac 12 \frac 12 1] &=& 
 \left(\ar{cccc} 1&0&0&0\\  0&\sqrt{\la}\,\sqrt{\frac{[p-1]}{[p]}} & 
 \sqrt{\la}\,\sqrt{\frac{[p+1]}{[p]}}&0\\ 0&0&0&1 \er\right) \ .\er  
\ee 

Now, straightforward check shows that the matrix $\F^{\frac 12 \frac 12}(p)$
found before can be obtained as follows:
$\F^{\frac 12 \frac 12}(p)=
{\rm C}'[\frac 12 \frac 12 0]{\rm C}_p[\frac 12\frac 12 0]+
{\rm C}'[\frac 12 \frac 12 1]{\rm C}_p[\frac 12\frac 12 1]$.
This suggests to consider the more general formula:
\be{7.4}
 \F^{IJ}(\p) \,=\, \sum_K {\rm C}'[IJK]\,{\rm C}_{\!\vec{\,p}}[IJK] \ .
\ee

Eq. (\ref{7.4}) may be regarded as an alternative definition of
the twisting element and it has already been considered in \cite{BS}
(similar expressions also appeared in \cite{BBB}) and proven to 
obey all axioms for the twisting element provided 
${\rm C}[IJK]$ and ${\rm C}_{\!\vec{\,p}}[IJK]$ are properly defined.
In this approach, however, the entries of the matrices
${\rm C}[IJK]$ and ${\rm C}_{\!\vec{\,p}}[IJK]$ are supposed to be
a-priory identified with some specific values of the CG coefficients 
and the $6j$-symbols (explicit formulae for them might be quite 
cumbersome). Moreover, these values (in general defined not uniquely) 
have to be chosen, to be compatible, in particular, with the choice of 
the matrices $R_\pm$ and $\R_\pm(\p)$. 
All this explains why we had not chosen eq. (\ref{7.4}) as a starting
point for constructing $\F^{\frac 12 \frac 12}(p)$. 

To sum up, we have demonstrated that eq. (\ref{7.4}) with ${\rm C}[IJK]$ 
and ${\rm C}_{\!\vec{\,p}}[IJK]$, built up from eigenvectors of the 
projectors $P^{IJ}_K$ and $\P^{IJ}_K$ according to given above 
prescriptions, gives correct expression for the twisting element. 
With this clarification the practical application of (\ref{7.4}) becomes
more straightforward.

Additionally, our approach allows to make the algebraic sense of 
(\ref{7.4}) more transparent. Indeed, since $P^{IJ}_K=$ 
${\rm C}'[IJK]{\rm C}[IJK]$, we can rewrite the formula for decomposition 
of $R$-matrices over the projectors (we used its simplest case (\ref{6.4}) 
above) in the  following form:
\be{7.5}
 \wh{R}^{IJ}_\pm \,=\, \sum_K {\rm C}'[IJK] \, r^{IJ}_{K,\pm} 
 \,{\rm C}[IJK] \ ,
\ee      
where $r^{IJ}_{K,\pm}$ are the corresponding eigenvalues 
[see \cite{FRT} for the fundamental representations of $U_q(sl(n))$ 
and  \cite{KS} for the highest irreps of $U_q(sl(2))$]. 
Now, bearing in mind the properties (\ref{2.8}) of the CG maps and
employing (\ref{7.4}), we can transform, according to (\ref{4.7}), eq.
(\ref{7.5}) into similar one for the twisted $R$-matrices:
\be{7.6}
 \wh{\R}^{IJ}_\pm(\p) \,=\, 
 \sum_K {\rm C}_{\!\vec{\,p}}'[IJK] \, r^{IJ}_{K,\pm} 
 \,{\rm C}_{\!\vec{\,p}}[IJK] \ .
\ee    
Thus, being written in the language of the projectors, the objects 
belonging to Hopf and quasi-Hopf structures look quite identical.

\subsection{On properties of the element $\chi$}
In conclusion, we wish to discuss in more detail properties of
the element $\chi$ which, as we have seen in \S 2.3 and \S 3.1, plays
an essential role in the theory of exact generating matrices.
First, we infer from (\ref{*D}) that the self-conjugated element 
{$\chi=\F(\p)\F^*(\p)$} obeys the relation 
\be{cD}
 \chi \, \D' \, = \, \D \, \chi \ . 
\ee
Additionally, eq. (\ref{ch}) also implies that $[R_\pm , \chi \chi' ]=0$. 
Together with (\ref{cD}) this allows us to assume that
\be{chi'}
   \chi' \,=\, \chi^{-1} \ .
\ee
Indeed, as has been demonstrated in \cite{BS}, there exists the
following universal expression for the element $\chi\in\G\o\G$:
\be{unc}
   \chi \,=\, \D(\k^{-1})\,(\k\o\k)\,R_+^{-1} \,=\,
    \D(\k)\,(\k\o\k)^{-1}\,R_-^{-1} \ ,
\ee
where $\k^2=v$ with $v$ being a certain invertible central element of $\G$
(the {\em ribbon element}, see \cite{rib}), such that
\be{rib}
  R_-^{-1}\,R_+ \,=\, \D(v^{-1})\,(v\o v) \ , \ \ \ 
  S(v) \,=\, v \ , \ \ \ \eps(v) \,=\, 1, \ \ \ v^* \,=\, v^{-1} \ .
\ee  
It is interesting to mention that, since the first relation in eq. 
(\ref{rib}) can be rewritten as 
\hbox{$\D(v^{-1})\,(v\o v)=\wh{R}_+^2=\wh{R}_-^{-2}$}
[recall that $\wh{R}_\pm \equiv PR_\pm$], eq. (\ref{unc}) admits the
following form:
$$
 \wh{R}_+ \,=\, \wh{\chi}^{-1} \, \bigl\{ \wh{R}_+^2 \bigr\}^{\frac 12} 
 \ , \ \ \
 \wh{R}_- \,=\, \bigl\{ \wh{R}_-^2 \bigr\}^{\frac 12} \, \wh{\chi} \ ,
$$ 
with $\wh{\chi}\equiv \chi P$. In other words, $\wh{\chi}^{\mp 1}$ appears 
as a (matrix) phase which fixes the choice of square root of 
$\wh{R}_\pm^2$. And, although, $\wh{\chi}^2=e\o e$ 
[according to (\ref{chi'})], this phase
turns out to be quite nontrivial, as we shall see below.

The properties (\ref{ch}) and (\ref{cD})-(\ref{chi'}) [as well as 
$\chi^*=\chi$]
become quite obvious for $\chi$ being defined as in (\ref{unc}).  
However, eq. (\ref{unc}) is not convenient if we want to get an explicit
form of $\chi^{IJ}$. Therefore, below we shall discuss, exploiting the 
language of quantum projectors, an alternative way of constructing $\chi$.

First, let us find with the help (\ref{6.11}) and (\ref{asF})
an explicit expression for the element $\chi$ in the fundamental
representation of $U_q(sl(2))$:
\be{6.12}
  \chi^{\frac 12 \frac 12} \,=\, 
 \left( \ar{cccc} 1&&& \\ [-1mm] & 2\la & -\om\la & \\ [-0.5mm]
  & \om\la & \ \ 2\la & \\ [-1mm] &&& 1 \er \right) \ .
\ee
In the non-deformed limit, $q\rar 1$, we have 
$\chi^{\frac 12 \frac 12} \rar e\o e$, 
as expected. Now we  notice that (\ref{6.12}) looks
very simply in terms of projectors (\ref{6.5}). Namely, for 
$\wh{\chi}^{\frac 12 \frac 12}=\chi^{\frac 12 \frac 12}
 P^{\frac 12 \frac 12}$ we have 
\be{ch'}
 \wh{\chi}^{\frac 12 \frac 12} \,=\, P_+ - P_- \ ,
\ee 
where $P_+\equiv P^{\frac 12 \frac 12}_{\ \,1}$, 
$P_-\equiv P^{\frac 12 \frac 12}_{\ \,0}$.
To explain this formula, we need the relation (\ref{7.4}). 
Substituting the latter into the definition of $\chi$, we obtain:
\be{cP}
  \wh{\chi}^{IJ} \,=\, \F^{IJ}\,(\F^{IJ})^* \, P^{IJ}\,=\, 
 \sum\nolimits_K \, C'[IJK] \, \overline{C[JIK]} \ ,
\ee
where $\overline{C[\dots]}=$ denotes a matrix complex conjugated to 
$C[\dots]$ and we have used the identity 
${\rm C}_{\!\vec{\,p}}[IJK]({\rm C}_{\!\vec{\,p}}[IJL])^*=\d_{KL}$
which follows from (\ref{*D}). Now, taking into account the well-known
symmetry of $q$-deformed CG coefficients (see, e.g., \cite{Bie})
\be{qq}
 C[JIK]_{q^{-1}}  \,=\, (-1)^{I+J-K}\,C[IJK]_q \ ,
\ee
we rewrite (\ref{cP}) as follows 
\be{FP}
 \wh{\chi}^{IJ} \,=\, 
 \sum\nolimits_{K=|I-J|}^{I+J} \, (-1)^{I+J-K}\, P^{IJ}_K \ .
\ee
Thus, for $U_q(sl(2))$ the element $\wh{\chi}$ is an altered sum of
the quantum projectors. For $U_q(sl(n))$ the symmetries of the CG
coefficients are of more sophisticated form (see, e.g., \cite{Kli}).
Therefore, in general, we should expect more complicated formula for 
$\chi$, but, presumably, still in terms of the quantum projectors.
In particular, formula (\ref{ch'}) remains, probably, true for any 
fundamental representation of $U_q(sl(n))$ [since in that case 
there are only two projectors and the corresponding coefficients
are uniquely fixed by the properties of $\chi$].

\vspace*{3mm}

{\bf Conclusion }
\vspace*{2mm}

 In the present paper we have demonstrated that the theory of (deformed) 
tensor operators and, in particular, the fusion procedure can be most 
naturally described employing the $R$-matrix language and revealing the 
underlying quasi-Hopf-algebraic structure. We clarified the role in
this context of the projectors and their $p$-dependent counterparts
which appear, respectively, in decompositions of $R$-matrices and 
twisted $R$-matrices.  From the practical point of view, the suggested 
prescription for constructing exact generating matrices 
can be used, e.g., for explicit computations and studies
of (deformed) CG coefficients for quantum Lie algebras of higher ranks.
On the other hand, the specific quasi-Hopf algebra [defined by the pair
$R_\pm$ and $\F(\p)$] appearing in this context
should certainly be studied in more detail since it provides
non-trivial (and presumably somewhat simplified) realization of the
abstract general scheme. Explicit formulae like that we have derived for 
$\F^{\frac 12 \frac 12}(p)$ may be useful here.

Although the present paper dealt mainly with the mathematical side of
the theory of tensor operators, we are going to discuss some physical
applications in future. 
{}Finally, we would like to note that it would be interesting to
extend the developed technique to the case of $q$ being a root of unity,
which would involve truncated quasi-Hopf algebras. 
\vspace*{2mm}

{\large\bf Acknowledgments }
\vspace*{2mm}

I am grateful to L.D.Faddeev, P.P.Kulish and V.Schomerus 
for useful discussions.
I am thankful to Prof.~R.Schrader for hospitality at 
Institut f\"ur Theoretische Physik, Freie Universit\"at Berlin
and to the T\"opfer Stiftung for financial support.

\newcommand{\sbibitem}[1]{\vspace*{-1.5ex} \bibitem{#1}}
 
\end{document}